\documentclass[letterpaper, 10 pt, conference]{ieeeconf}  
\IEEEoverridecommandlockouts                              

\usepackage{graphics} 
\usepackage{epsfig}   
\usepackage{mathptmx} 
\usepackage{times}    
\usepackage{amsmath}  
\usepackage{amssymb}  
\usepackage{graphicx}
\usepackage{svg}
\usepackage{caption}
\usepackage{subcaption}
\usepackage{cite}
\usepackage{xcolor}
\usepackage{tcolorbox}
\usepackage{algorithm}
\usepackage{algpseudocode}
\usepackage[super]{nth}
\usepackage[utf8]{inputenc}
\usepackage{upgreek}
\usepackage[scr]{rsfso}
\usepackage{bm}       
\usepackage{dsfont}   
\usepackage{float}
\usepackage{soul}
\usepackage[percent]{overpic}
\usepackage{anyfontsize}
\setstcolor{red}

\makeatletter
\let\NAT@parse\undefined
\makeatother
\usepackage[colorlinks=true, allcolors=blue]{hyperref}

\usepackage[siunitx]{circuitikz}
\usepackage{tikz}
\usetikzlibrary{decorations.shapes}
\tikzset{decorate sep/.style 2 args=
{decorate,decoration={shape backgrounds,shape=circle,shape size=#1,shape sep=#2}}}

\definecolor{tueRed}{HTML}{C81919}
\definecolor{mlBlue}{HTML}{0072BD}
\definecolor{mlOrange}{HTML}{D95319}

\title{\LARGE \bf
Constraints-Informed Neural-Laguerre Approximation of Nonlinear MPC with Application in Power Electronics
}

\author{Duo Xu$^{1}$, Rody Aerts$^{1}$, Petros Karamanakos$^{2}$, Mircea Lazar$^{1}$
\thanks{This work is funded by the EU Horizon 2020 research project IT2 (IC Technology for the 2nm Node), Grant agreement ID: 875999.}
\thanks{$^{1}$D.~Xu, R.~Aerts and M.~Lazar are with the Department of Electrical Engineering, Eindhoven University of Technology, P.O. Box 513, 5600 MB Eindhoven, The Netherlands
{\tt\small d.xu@tue.nl; r.j.c.aerts@student.tue.nl; m.lazar@tue.nl}}%
\thanks{$^{2}$P.~Karamanakos is with the Faculty of Information Technology and Communication Sciences, Tampere University, 33101 Tampere, Finland
{\tt\small p.karamanakos@ieee.org}}%
}

\begin{document}

\maketitle
\thispagestyle{empty}
\pagestyle{empty}

\begin{abstract}
This paper considers learning online (implicit) nonlinear model predictive control (MPC) laws using neural networks and Laguerre functions. Firstly, we parameterize the control sequence of nonlinear MPC using Laguerre functions, which typically yields a smoother control law compared to the original nonlinear MPC law. Secondly, we employ neural networks to learn the coefficients of the Laguerre nonlinear MPC solution, which comes with several benefits, namely the dimension of the learning space is dictated by the number of Laguerre functions and the complete predicted input sequence can be used to learn the coefficients. To mitigate constraints violation for neural approximations of nonlinear MPC, we develop a constraints-informed loss function that penalizes the violation of polytopic state constraints during learning. Box input constraints are handled by using a clamp function in the output layer of the neural network. We demonstrate the effectiveness of the developed framework on a nonlinear buck-boost converter model with sampling rates in the sub-millisecond range, where online nonlinear MPC would not be able to run in real time. The developed constraints-informed neural-Laguerre approximation yields similar performance with long-horizon online nonlinear MPC, but with execution times of a few microseconds, as validated on a field-programmable gate array (FPGA) platform.
\end{abstract}
\section{Introduction}

In contemporary power electronics, model predictive control (MPC) has gained substantial attention and application. 
The imperative for MPC lies in computing actions within sub-microsecond intervals to meet the stringent requirements of operation at high switching frequencies. The nonlinear dynamics inherent in power electronic systems, however, pose challenges from a control perspective, adding further computational complexities \cite{karamanakos2020model}. 

To address the computational efficiency problem, neural networks (NN) are proposed to approximate the control law of MPC, which was first studied in 1995 \cite{parisini1995receding}.
It has been shown that NNs can exactly represent the explicit MPC law
\cite{karg2020efficient} and the feasibility can be enhanced by using Dykstra's projection algorithm \cite{chen2018approximating}.
In \cite{maddalena2020neural}, a parametric quadratic program layer is introduced into the NN to simplify the explicit MPC. In \cite{drgona2020learning}, a differentiable MPC framework was proposed as an unsupervised, sampling-based approximate solution. This framework addresses the MPC problem by considering state and input constraints as soft constraints, with violations accounted for within the loss function formulation.
For nonlinear systems, the explicit nonlinear MPC (NMPC) solution is not straightforward, which poses extra challenges in approximation via NN. 
In \cite{hertneck2018learning}, the NMPC is robustly designed with a chosen bound on the input error. If the approximation error of the NN control law is below the defined bound, the stability and constraint satisfaction are guaranteed. 
In \cite{drgovna2022differentiable}, the differentiable MPC problem was extended to accommodate nonlinear systems, which are identified based on neural state-space model architecture.

Among most of the works dealing with MPC approximation, the first input action of the MPC is learnt instead of the full predicted input sequence. This is intuitive because only the first input action is applied to regulate the system.
In recent works, a safety-augmented NN is proposed to approximate the full predicted input sequence of NMPC, which achieves closed-loop control performance with guaranteed safety \cite{hose2023approximate}. 
Moreover, the approximated full input sequence, considered a sub-optimal solution, contributes to ensuring asymptotic stability with high probability \cite{wang2023policy}.
Thus, the inclusion of safety and stability concerns underscores the significance of approximating the full input sequence.

In \cite{wang2004discrete}, an orthonormal function (Laguerre function) was proposed to parametrize the entire predicted input sequence of the linear MPC problem. The input sequence is approximated through a weighted sum of a reduced number of Laguerre coefficients. The Laguerre coefficients are solved instead with fewer numbers of variables, while the full input sequence is stored inherently.
In \cite{bamimore2023laguerre}, the input sequence parametrization using Laguerre function was extended to a nonlinear model.
As observed, the Laguerre function parameterized NMPC (LagNMPC) results in smoother and more cautious input sequences compared to the original NMPC.

Considering the above, this paper explores the approximation of LagNMPC using NNs, allowing for the learning of the entire predicted control sequence. An additional objective is to ensure that the resulting approximate explicit NN controllers satisfy input and state constraints. The main contributions are: 
\emph{(i)}~an NN architecture for learning online Laguerre NMPC control sequences; 
\emph{(ii)}~for power electronic converters, hard input constraints guarantees via clamp function in the output layer.
Similar to the soft constraints used in the unsupervised approximation scheme proposed in \cite{drgona2020learning}, a constraints-informed loss function is employed during learning, which penalizes violations of the polytopic constraints.
Furthermore, to mitigate the steady-state error due to the learning error, we adopt an offset-free based NN control law inspired by \cite{chan2021deep}. To demonstrate the effectiveness of the proposed method, we apply it to a case study from the field of power electronics, namely a buck-boost dc-dc converter, with the sampling interval of $100$ microseconds and for a prediction horizon of $20$ time steps. As shown, the developed neural explicit approximations of LagNMPC can run on a field-programmable gate array (FPGA) in under $20$ microseconds without parallelization, and under $3$ microseconds with parallelization. 


\paragraph*{Notation and basic definitions} Let the symbols $\mathbb{R}$, $\mathbb{Z}$ and $\mathbb{Z}_+$ denote the set of real, integer and non-negative integer numbers, respectively. 
$\mathbb{Z}_{[i,j]}$ denotes the set of integer numbers restricted in $[i,j]$ ($i < j$).
$\Vert x \Vert^2_P$ denotes the quadratic form $x^TPx$.
$\lceil \cdot \rceil$ denotes the ceil function.
The function $f$ of $g(x)$ will be denoted as $f\circ g(x) = f(g(x))$.
\section{Nonlinear Model Predictive Control and Laguerre Function Parameterization} \label{Section2}

\subsection{Nonlinear model predictive control}
In this paper, the following discrete-time nonlinear system is considered:
\begin{equation} \label{eq2:model}
    x_{k+1} = f(x_k, u_k),\quad k\in\mathbb{Z}_+,
\end{equation}
where $x_k \in \mathbb{R}^{n_x}$ is the state and $u_k \in \mathbb{R}^{n_u}$ the control inputs. The steady-state solution of the system exists with $x^{\text{ss}} = f(x^{\text{ss}}, u^{\text{ss}})$.
The NMPC optimization problem is formulated as
\begin{subequations} \label{eq2:NMPC}
    \begin{align}
    \min_{U_{k}}\;\; &  J(x_k,U_k) = ||x_{N|k}-x^{\text{ss}}||_P^2 \nonumber\\ 
    &+\sum_{i=0}^{N-1} ||x_{i|k}-x^{\text{ss}}||_Q^2 + ||u_{i|k}-u^{\text{ss}}||_R^2  \\
    \text{s.t. }& x_{i+1|k} = f(x_{i|k},u_{i|k}), &&\hspace{-5em} \forall i \in\mathbb{Z}_{[0,N-1]},\\ 
    & x_{i|k} \in \mathbb{X}, &&\hspace{-5em} \forall i \in \mathbb{Z}_{[1,N]}, \\
    & u_{i|k} \in \mathbb{U}, &&\hspace{-5em} \forall i \in \mathbb{Z}_{[0,N-1]},
    \end{align}
\end{subequations}
with $x_{0|k} = x_k $ the measured state, $U_{k}=[u_{0|k}^{\top} , \ldots,  u_{N-1|k}^{\top} ]^{\top}$ the predicted input sequence over a prediction horizon $N$, $\{P,Q \in \mathbb{R}^{n_x \times n_x}\;|\; P,Q \succeq 0\}$ the penalty on states and $\{R \in \mathbb{R}^{n_u \times n_u}\;|\; R \succ 0\}$ the penalty on inputs. The state and input constraint sets are represented by $\mathbb{X}$ and $\mathbb{U}$ respectively. 
The optimal sequence of the predicted input $U_{k}^*=[ u_{0|k}^{*{\top}}, \ldots, u_{N-1|k}^{*{\top}}]^{\top}$ obtained after solving the optimization problem \eqref{eq2:NMPC} yields the first element $u_{0|k}^*$ as the optimal input. In general, the control law of the NMPC is implicitly given as
\begin{equation} \label{eq2:u_NMPC}
    u_k^* = u_{0|k}^{*} = \pi(x_k).
\end{equation}
A set of terminal sets and terminal cost functions can be derived, such that the stability and recursive feasibility of the NMPC problem \eqref{eq2:NMPC} can be guaranteed \cite{eyubouglu2022snmpc}.

\subsection{Parameterized NMPC using Laguerre function}
In \cite{wang2004discrete}, a discrete-time Laguerre function is introduced to reformulate the NMPC problem and simplify the solution. In such way, a long prediction horizon can be achieved without using a large number of parameters. Moreover, Laguerre function is known to have rapid exponential decay \cite{bamimore2023laguerre}. In this section, a single-input system ($n_u=1$) parameterization using the Laguerre function is considered. Nevertheless, the extension to multi-input ($n_u>1$) formulation is straightforward, as demonstrated in \cite{wang2009model}.

Consider a sequence of the predicted inputs $U_k$ of the NMPC problem and the steady-state input vector $U^{\text{ss}} = [ u^{\text{ss}}  , \ldots , u^{\text{ss}} ]^\top$ with the same dimension $N$. By introducing the Laguerre function, the error between $U_k$ and $U^{\text{ss}}$ can be parameterized as
\begin{equation}
\begin{aligned}
U_k-U^{\text{ss}} &  = L \eta_k  \\
&= 
\begin{bmatrix} 
l_1(0) &  \ldots & l_M(0)  \\
\vdots &        & \vdots \\
l_1(N-1) &  \ldots & l_M(N-1) \\
\end{bmatrix}
\begin{bmatrix} \eta_{k,1} \\ \vdots \\ \eta_{k,M} \end{bmatrix},
\end{aligned}
\end{equation} 
where $\eta_k = \begin{bmatrix} \eta_{k,1} & \ldots & \eta_{k,M}\end{bmatrix}^\top \in \mathbb{R}^M$ and $M$ is the enumeration of the orthogonal functions $L\in \mathbb{R}^{N\times M}$. 
Each column of $L$ is defined as the orthonormal basis function $\boldsymbol{l}_\rho$ given by
\begin{equation}
\boldsymbol{l}_\rho = \begin{bmatrix}
    l_\rho(0) & l_\rho(1) & \ldots & l_\rho(N-1)
\end{bmatrix}^\top, \forall \rho\in \mathbb{I}_{[1,M]}.
\end{equation} 
The predicted input sequence $U_k$ is posed as a weighted sum of a certain number of orthonormal basis functions as
\begin{equation}
    U_k-U^{ss} = \sum_{\rho=1}^M \boldsymbol{l}_\rho \eta_{k,\rho}.
\end{equation} 
In addition, by defining the transpose of each row of $L$ as
\begin{equation}
    L_i = \begin{bmatrix}
    l_1(i) & l_2(i) & \ldots & l_M(i)
\end{bmatrix}^\top,
\end{equation}
the Laguerre function can be expressed by the following difference equation
\begin{equation}
    L_{i+1} = A_L L_i,
\end{equation} 
\begin{equation}
     A_L = \begin{bmatrix}
    \alpha  \\
    \beta & \alpha \\
    -\alpha \beta & \beta & \alpha \\
    \vdots & \vdots & \ddots& \ddots \\
    (-\alpha)^{M-2}\beta & (-\alpha)^{M-3}\beta & \ldots & \beta & \alpha
    \end{bmatrix},
\end{equation} 
where $\alpha\in[0,1)$ is the pole of the discrete-time Laguerre network, $\beta=(1-\alpha^2)$ and
$L_0=\sqrt{\beta} \begin{bmatrix} 1 & -\alpha & \alpha^2 & \ldots & (-\alpha)^{M-1}\end{bmatrix}^\top$. Therefore, the parameterized predicted input sequence for the NMPC problem is derived as
\begin{equation}
U_k = \begin{bmatrix}u_{0 \mid k} \\ u_{1 \mid k} \\ \vdots \\ u_{N-1 \mid k} \end{bmatrix} = 
\begin{bmatrix} L_0^\top \eta_k + u^{\text{ss}} \\ L_1^\top \eta_k + u^{\text{ss}}\\ \vdots \\ L_{N-1}^\top \eta_k + u^{\text{ss}}\end{bmatrix}.
\end{equation} 

The LagNMPC problem can be formulated as
\begin{subequations} \label{eq3:LagNMPC}
    \begin{align}
    \min_{\eta}\;\; &  \Vert x_{N|k}-x^{\text{ss}} \Vert^2_P  +\sum_{i=0}^{N-1} \Vert x_{i|k}-x^{\text{ss}}\Vert^2_Q + \Vert L_i^\top \eta_k \Vert^2_R \\
    \text { s.t. } 
    & x_{i+1|k}=f(x_{i|k},L_i^\top \eta_k+u^{\text{ss}}),  &&\hspace{-5em}\forall i \in \mathbb{Z}_{[0,N-1]}, \label{equ2:prediction1}\\
    & L_i^\top \eta_k +u^{\text{ss}}\in \mathbb{U},  &&\hspace{-5em}\forall i \in \mathbb{Z}_{[0,N-1]}, \\
    & x_{i|k} \in \mathbb{X},  &&\hspace{-5em}\forall i \in \mathbb{Z}_{[1,N]}.
    \end{align}
\end{subequations}
If the Laguerre function is set with $M=N$ and $\alpha=0$, the LagNMPC problem is equivalent with the original NMPC problem \cite{wang2009model}.
Consequently, the optimal input can be recovered by
\begin{subequations} \label{eq2:U_LagNMPC}
    \begin{align}
    U_k^* &= L \eta_k^* + U^{\text{ss}},\label{eq3:Ukhat}\\
    u_k^* &= u_{0|k}^* = L_0^\top \eta_k^* + u^{\text{ss}}. \label{eq2:u_LagNMPC}
    \end{align}
\end{subequations}
Similar to the NMPC problem, the LagNMPC optimization problem \eqref{eq3:LagNMPC} implicitly introduces the following solution
\begin{equation} \label{eq2:eta_LagNMPC}
    \eta_k^*  = \Pi_\eta(x_k).
\end{equation}
Still, a set of terminal sets and terminal cost functions can be derived, such that the stability and recursive feasibility of the LagNMPC problem \eqref{eq2:NMPC} can be guaranteed \cite{eyubouglu2022snmpc}.
\section{Approximate NMPC using Neural Networks with Hard Input Constraints Guarantees} \label{Section3}

\subsection{Data generation}
To approximate either the NMPC or the LagNMPC control policy, data need to be collected to train the NNs. Thus, we propose the quasi-random Halton sequence to sample the constraint set of states for data collection. The Halton sequence spanning $d$ dimensions is defined as \cite{cheng2013computational}
\begin{equation}
    H(d,n) = \big( \phi_{p_1}(n),...\;,\phi_{p_j}(n),...\;,\phi_{p_d}(n) \big),
    \label{eq3:xn}
\end{equation}
with $p_1,...\;,p_j,...\;,p_d$ being the first $d$ prime numbers, $n \in \mathbb{Z}$ and $\phi_{p_j}(n) \; : \; \mathbb{Z} \rightarrow [0,1)$ the radical inverse function
\begin{equation}
    \phi_{p_j}(n) = \sum^{l(j)}_{i=0} \frac{a_i(j,n)}{p_j^{i+1}}, \hspace{3mm} l(j) = \lceil \text{log}_{p_j}(n) \rceil,
    \label{eq3:phi_pj}
\end{equation}
where $a_j(j,n) \in [0,p_j-1]$ the integer coefficients from the $p_j$-ary expansion of
\begin{equation}
    n = \sum^{l(j)}_{i=0} a_i(j,n) p_j^i.
    \label{eq3:expa}
\end{equation}

Samples obtained with (\ref{eq3:xn}) are located in the $d$-dimensional hypercube $[0,1)^d$. To transform the samples to $\mathbb{X}$ a bounding hypercube $\mathbb{B}$ is obtained via
\begin{equation}
    \min_{\mathbb{B}} |\mathbb{B} - \mathbb{X} | \text{ s.t. } \mathbb{X}  \subseteq \mathbb{B}.
    \label{eq3:boundingbox}
\end{equation}
Hence, with the bounds of $\mathbb{B}$ known, the samples are scaled and translated to match $\mathbb{X}$. 

Algorithm \ref{alg:Halton} describes how the states are sampled. It uses the constraint set of states $\mathbb{X}$ and the number of desirable data points $N_d$. First the hypercube $\mathbb{B}$ bounding $\mathbb{X}$ is obtained, next samples are obtained using the Halton sequence for the stated set of integers which are scaled and translated. Only the states within $\mathbb{X}$ are added to the dataset and the process is repeated until $N_d$ samples are collected, creating a set of uniformly distributed states. 

\begin{algorithm}[t!]
    \caption{Halton sampling algorithm.}
    \begin{algorithmic}[1]   
    \Function{Halton}{$\mathbb{X}$, $N_d$}
    \State obtain $\mathbb{B}$
    \State $n_d,j \gets 0$
    \State $s \gets 1$
    \While{$n_d < N_d$}
    \State $\tilde{x} \gets H(n_x,\mathbb{Z}_{[N_d j, \;N_d(j+1)-n_d]})$
    \State scale and translate $\tilde{x}$ 
    \State $n_d \gets n_d + \big|\{\tilde{x}$ $|$ $\tilde{x} \in \mathbb{X}\}\big|$
    \State $\{x_i\}_{i=s}^{n_d} \gets \{\tilde{x}$ $|$ $\tilde{x} \in \mathbb{X}\}$
    \State $s \gets n_d + 1$
    \State $j \gets j + 1$    
    \EndWhile
    \State \Return $x$
    \EndFunction
    \end{algorithmic} 
    \label{alg:Halton}
\end{algorithm}
After a valid set of states $x_i$ is sampled using the Halton method within the state constraint set $\mathbb{X}$,
the corresponding optimal solution $\eta_i^*$ ($u_i^*$ for NMPC) is obtained by solving $ \eta_i^* = \Pi_\eta(x_i) $ ($u_i^* =\pi_(x_i)$ for NMPC).
A total number of $N_s$ ($N_s\leq N_d$) dataset $\{x_i, \eta_i^*\}_{i=1}^{N_{s}}$ (or $\{x_i, u_i^*\}_{i=1}^{N_{s}}$) is collected when the infeasible solutions of the (Lag)NMPC problem are removed.

\subsection{Neural network architectures}
\subsubsection{NMPC}

\begin{figure}[t!]
    \centering
    \includegraphics[width=\linewidth]{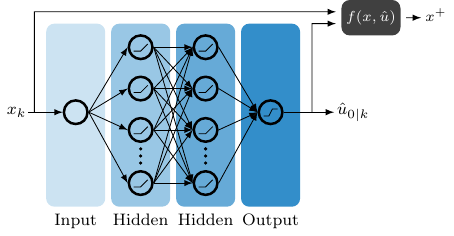}
    \caption{MLP diagram of NN-based NMPC: $\hat{u}=\hat{\pi}(x,\theta)$.}
    \label{fig2:NN}
\end{figure}

The multi-layer perceptron (MLP) neural network is adopted to approximate the NMPC solution as depicted in Figure \ref{fig2:NN}.
The MLP network is a feedforward NN complement with $h$ hidden layers and $n$ nodes in each hidden layer, and the nodes of consecutive layers are fully connected. 
The proposed MLP network learns a parameterized function $\hat{u} = \hat{\pi}(x,\theta) $ and the approximated control action is recovered by
\begin{equation} \label{eq2:MLP}
\begin{aligned}
\hat{u} = \hat{\pi}(x,\theta)  = r_c \circ g_{h+1} \circ r \circ g_{h} \circ \cdot\cdot\cdot \; \circ r \circ g_1(x),    
\end{aligned}
\end{equation}
where $g_{h}(\xi)= W_{h} \xi + b_{h}$ is the affine function, $r(\xi)=\max\{0,\xi\}$ the element-wise rectified linear units (ReLU) activation function, and
$\theta=\{W_j,b_j\}_{j=1}^{h+1}$ represents the parameters to be optimized. 
Considering box constraints on the input, the clamp activation function is applied element-wise to the output of the NN and is defined as follows:
\begin{equation} \label{eq3:activation_clamp}
\begin{aligned}
r_c(\xi,u_{\min},u_{\max}) =& \min{\{\max{\{\xi,u_{\min}\}},u_{\max}\}} \\
 =& \begin{cases} \begin{aligned}
u_{\min} \quad& \text{if} \quad \xi<u_{\min}, \\
\xi \quad& \text{if} \quad u_{\min}\leq\xi\leq u_{\max}, \\
u_{\max} \quad& \text{if} \quad \xi>u_{\max}.
\end{aligned} \end{cases}
\end{aligned}
\end{equation}
This function is bounded within the interval $[u_{\text{min}}, u_{\text{max}}]$, where $u_{\text{min}}$ and $u_{\text{max}}$ represent the lower and upper bounds on the input. As a result, the box input constraints are Incorporated during the training process of the network.

In order To train the MLP network, the following loss function is formulated:
\begin{equation} \label{eq3:lossNMPC}
    \mathscr{L}_{s}(\theta) = 
    \frac{1}{N_s} \sum_{i=1}^{N_s} \Vert u_i^* - \hat{\pi}(x_i,\theta) \Vert_2^2.
\end{equation}
This loss function considers the mean-squared error between the optimal input sequence $u_i^*$ and the output of the network $\hat{\pi}(x_i,\theta)$. 

Once the training solution $\hat{\theta}$ is derived, 
the approximated NN-based NMPC law $\hat{u}=\hat{\pi}(x,\hat{\theta})$ is implemented instead, which gives close control performance with greatly reduced computational time.

\subsubsection{LagNMPC}

\begin{figure}[t!]
    \centering
    \includegraphics[width=\linewidth]{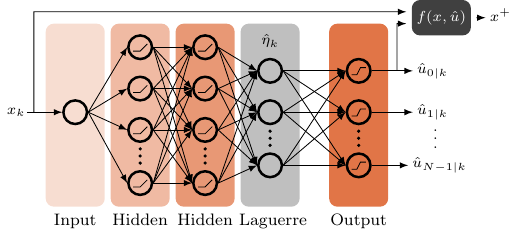}
    \caption{MLP diagram of NN-based LagNMPC: $\hat{U} = \hat{\Pi}(x,\theta,L,U^{\text{ss}})$.}
    \label{fig:NN_Lag}
\end{figure}

In order to approximate LagNMPC solution \eqref{eq2:eta_LagNMPC} using NN, we modify the MLP network structure, which is visualized in Figure \ref{fig:NN_Lag}. The Laguerre layer represents the operation of \eqref{eq3:Ukhat} and can be seen as an affine operation of the network with fixed weights and biases. The clamp activation function \eqref{eq3:activation_clamp} is applied element-wise to the results of the Laguerre layer. The proposed MLP network learns a parameterized function $\hat{U} = \hat{\Pi}(x,\theta,L,U^{\text{ss}}) $, and the approximated control action sequence is recovered by
\begin{subequations}
    \begin{align}
    &\hat{U} = \hat{\Pi}(x,\theta,L,U^{\text{ss}}) = r_c \circ g_L \circ g_{h+1} \circ r \circ g_{h} \circ \cdot\cdot\cdot \; \circ r \circ g_1(x),\\
    &g_L(\hat{\eta}) = L\hat{\eta} + U^{\text{ss}},
    \end{align}
\end{subequations}
where $g_L(\hat{\eta})$ represents the Laguerre layer. If only the first entry of the input sequence $\hat{u}$ is needed for control, the proposed NN architecture can be simplified as follows:
\begin{subequations}
    \begin{align}
    &\hat{u} = \hat{\Pi}_0(x,\theta,L_0,u^{\text{ss}}) = r_c \circ g_{L_0} \circ g_{h+1} \circ r \circ g_{h} \circ \cdot\cdot\cdot \; \circ r \circ g_1(x),\\
    &g_{L_0}(\hat{\eta}) = L_0^\top\hat{\eta} + u^{\text{ss}},
    \end{align}
\end{subequations}
where the Laguerre layer can be reduced to $g_{L_0}(\hat{\eta})$.

Given that the optimal input sequence is recovered by $U^*= L \eta^*+U^{\text{ss}}$, the following loss function is proposed to minimize the estimation error in the predicted input vector space:
\begin{equation}
    \label{eq3:lossLagNMPC}
    \mathscr{L}_{s}(\theta) = \frac{1}{N_s} \sum_{i=1}^{N_s} \Vert L \eta^*_i +U^{\text{ss}} - \hat{\Pi}(x_i,\theta,L,U^{\text{ss}})\Vert^2_2.
\end{equation}

\subsection{Training of NN with state constraint information}
To improve the performance of the NN in approximating both NMPC and LagNMPC control laws and to reduce the risk of violating the state constraints, we introduce the state constraints-informed (ConInf) loss function as follows:
\begin{equation} 
    \mathscr{L}_x(\theta) = \frac{1}{N_s} \sum_{i=1}^{N_s}  \Vert \mathscr{L}_x^i(\theta) \Vert_{\Gamma}^2, \label{eq3:Loss_x}
\end{equation}
where $\Gamma$ is a scaling matrix. In this section, box constraints are imposed on the system state, which are expressed as
\begin{equation}
    \mathscr{L}_x^i(\theta) = \max\{ 0 , x_i^+(\theta) - x_{\max}  \}+ \min\{ 0 , x_i^+(\theta) - x_{\min}  \} ,
\end{equation}
where $x_i^+(\theta)=f( x_i , \hat{u}_i )$ represents the one-step predicted state resulting from the approximated input $\hat{u}_i$, which is defined as
\begin{equation} \label{eq3:NNcontroller}
\begin{aligned}
\hat{u}_i = 
\begin{cases}
    \hat{\pi}(x_i,\hat{\theta}) , & \text { if } \text{NMPC}, \\ 
    \hat{\Pi}_0(x_i,\hat{\theta},L_0,u^{\text{ss}}) , & \text { if } \text{LagNMPC}.
\end{cases}
\end{aligned}
\end{equation}
Note that a more general formulation for polytopic state constraints is available in \cite{drgona2020learning}. By adding \eqref{eq3:Loss_x} in the standard loss function, \eqref{eq3:lossNMPC} or \eqref{eq3:lossLagNMPC}, the constraints-informed (ConInf) loss function is obtained and used to train the MLP networks, i.e.,
\begin{equation} \label{eq3:ConInf}
    \mathscr{L}(\theta) = \mathscr{L}_s(\theta)  + \mathscr{L}_x(\theta).
\end{equation}
This formulation assists in respecting the state constraints by minimizing their violations, treating them  as soft constraints.
Additionally, incorporating the activation function \eqref{eq3:activation_clamp} into the NN output ensures that input constraints, considered as hard constraints, are strictly enforced.
 
Except for the constraints-informed losses, the error between the gradients of $u_i^*$ and $\hat{u}_i$ can also assist in improving the approximation performance \cite{Winqvist2021}. Nevertheless, the inclusion of such a loss function does not appear to be beneficial  for the problem considered in this work and is therefore omitted.

\subsection{Offset-free based control law}
As investigated, since there exists inevitable estimation error between the original (Lag)NMPC solution and the NN approximated solution, it is likely that
\begin{subequations} 
\begin{align}
\hat{\pi}(x^{\text{ss}},\hat{\theta})  &\neq u^{\text{ss}}, \\
\hat{\Pi}_0(x^{\text{ss}},\hat{\theta},L_0,u^{\text{ss}}) &\neq u^{\text{ss}}.
\end{align}
\end{subequations}
Thus, a steady-state solution with offset can be derived when the system is controlled via the NN-based controller. In order to handle this problem, an offset-free based control law is adopted inspired by \cite{chan2021deep}, i.e.,
\begin{equation} \label{eq3:offset_control}
\begin{split}
&\hat{u}_k^\text{of} = \\
& \begin{cases}
\pi(x_k,\hat{\theta}) - \pi(x^{\text{ss}},\hat{\theta}) + u^{\text{ss}}, &\hspace{-0.8em}\text{ if } \text{NMPC}, \\ 
\hat{\Pi}_0(x_k,\hat{\theta},L_0,u^{\text{ss}}) - \hat{\Pi}_0(x^{\text{ss}},\hat{\theta},L_0,u^{\text{ss}}) + u^{\text{ss}}, &\hspace{-0.8em}\text{ if }\text{LagNMPC}.
\end{cases}
\end{split}
\end{equation}
However, since the offset-free control law introduces a constant offset on the system input, the control action runs the risk of violating both state and input constraints. Therefore, the offset-free control law is applied when the state enters a set $\mathbb{X}_T=\{x_k\ :\ \Vert x_k - x^{ss}\Vert_2 \leq \epsilon\}$. Notice that typically the desires steady-state lies within the interior of the state constraints set $\mathbb{X}$.
\section{Application to a Buck-Boost Converter} \label{Section4}

\subsection{Buck-Boost DC-DC converter}
\begin{figure}[b!]
    \centering
    \includegraphics[width=\linewidth]{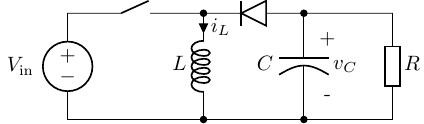}
    \caption{Schematic of the buck-boost converter.}
    \label{fig2:buckboost}
\end{figure}
\begin{table}[b!] 
\centering
\caption{Buck Boost Converter Parameters}
\label{tab2:Circuit Parameters}
\begin{tabular}{|l|l|l|}
\hline
Parameter & Symbol &  Value \\ \hline
Sampling interval & $T_s$ &  $100\upmu$s \\
Bus voltage & $V_{\text{in}}$ &  $15$V \\
Inductance & $L$ &  $4.2$mH \\
Capacitance & $C$ &  $2200\upmu$F \\
Load resistance & $R$ &  $165 \Omega$ \\
Constraint & $x_{\min}$ &  $\begin{bmatrix}0.01 & -20\end{bmatrix}^\top$ \\
Constraint & $x_{\max}$ &  $\begin{bmatrix}2 & 0\end{bmatrix}^\top$ \\
Constraint & $u_{\min}$ &  $0.1$ \\
Constraint & $u_{\max}$ &  $0.9$ \\
\hline
\end{tabular}
\end{table}
Figure \ref{fig2:buckboost} illustrates the schematic of an ideal dc-dc buck-boost converter, i.e., the parasitic components are neglected. The converter parameters are listed in Table \ref{tab2:Circuit Parameters}. 
The discrete-time nonlinear averaged model of the converter is of the form \cite{lazar2008input}
\begin{equation} \label{eq:buckboost}
\begin{aligned}
x_{k+1} = f(x_k,u_k) = \begin{bmatrix} 1 & \frac{T_s}{L} \\ -\frac{T_s}{C} & 1-\frac{T_s}{RC} \end{bmatrix} x_k + \begin{bmatrix} \frac{T_s}{L}(V_{\text{in}}-x_{2,k}) \\ \frac{T_s}{C} x_{1,k}\end{bmatrix} u_{k},
\end{aligned}
\end{equation}
where the state $x_k = \begin{bmatrix} x_{1,k} & x_{2,k} \end{bmatrix}^\top$ consists of the inductance current $x_{1,k} = i_L$, and the capacitor voltage $x_{2,k} = v_C$. Moreover, the duty cycle is treated as the control input $u_{k}$.


The control purpose is to regulate the capacitor voltage $x_{2,k}$ to a desired steady-state value $x_2^{\text{ss}}$. Accordingly, the steady-state value of $x_1^{\text{ss}}$ and $u^{\text{ss}}$ can be derived as
\begin{equation} \label{eq2:steadystate}
x_1^{\text{ss}} = \frac{x_2^{\text{ss}}}{R(u^{\text{ss}}-1)},\hspace{0.5cm} 
u^{\text{ss}} = \frac{x_2^{\text{ss}}}{x_2^{\text{ss}} - V_{\text{in}}}.
\end{equation}


\subsection{Controller and neural network setup}
The performance of both NMPC and LagNMPC and their corresponding NN-based controllers are simulated using MATLAB. The (Lag)NMPC problem is solved with the nonlinear solver \textsf{fmincon}. Below Table \ref{tab5:MPC_setup} lists the parameters we set for the two NMPC problems.
\begin{table}[h]
\centering
\caption{MPC \& LagNMPC setup}
\label{tab5:MPC_setup}
\begin{tabular}{|l|l|l|}
\hline
Parameter          & Symbol   & Value \\ \hline
Prediction horizon & $N$      & 20 \\
Weighting matrix   & $Q$      & diag(1, 0.1) \\
Weighting matrix   & $R$      & 0.7 \\
Weighting matrix   & $P$      & diag(10, 1) \\ 
Laguerre size      & $M$      & 4 \\
Laguerre pole      & $\alpha$ & 0.9 \\
Equilibrium point  & $x^{\text{ss}}$    & $\begin{bmatrix}0.101 & -10\end{bmatrix}^\top$ \\ 
Equilibrium point  & $u^{\text{ss}}$    & 0.4 \\ 
\hline
\end{tabular}
\end{table}

Algorithm \ref{alg:Halton} is implemented to uniformly sample the data of state inside the state constraint set $\mathbb{X}$. 
A total of 20000 pairs of the feasible dataset $\{x_i, u_i^*\}_{i=1}^{N_s}$ ($\{x_i, \eta_i^*\}_{i=1}^{N_s}$) are collected for approximating the (Lag)NMPC problem.

The networks were trained in PyTorch for 1000 epochs, with a learning rate of \num{1e-4}, the \textsf{AdamW} optimizer and using batch normalization with a batch size of 1024. Note that $70\%$ of the dataset is reserved for training the network, while the remaining is used for validation purpose. The weights in the loss function term (\ref{eq3:Loss_x}) are chosen to be $\Gamma=\text{diag}(1,0.1)$. The MLP network consists of $h=2$ hidden layers and $n=20$ nodes per hidden layer. 

\subsection{Performance comparison} 
\begin{figure*}[t]
    \centering
    \begin{subfigure}[b]{0.32\linewidth}
        \centering
        \includegraphics[width=\linewidth]{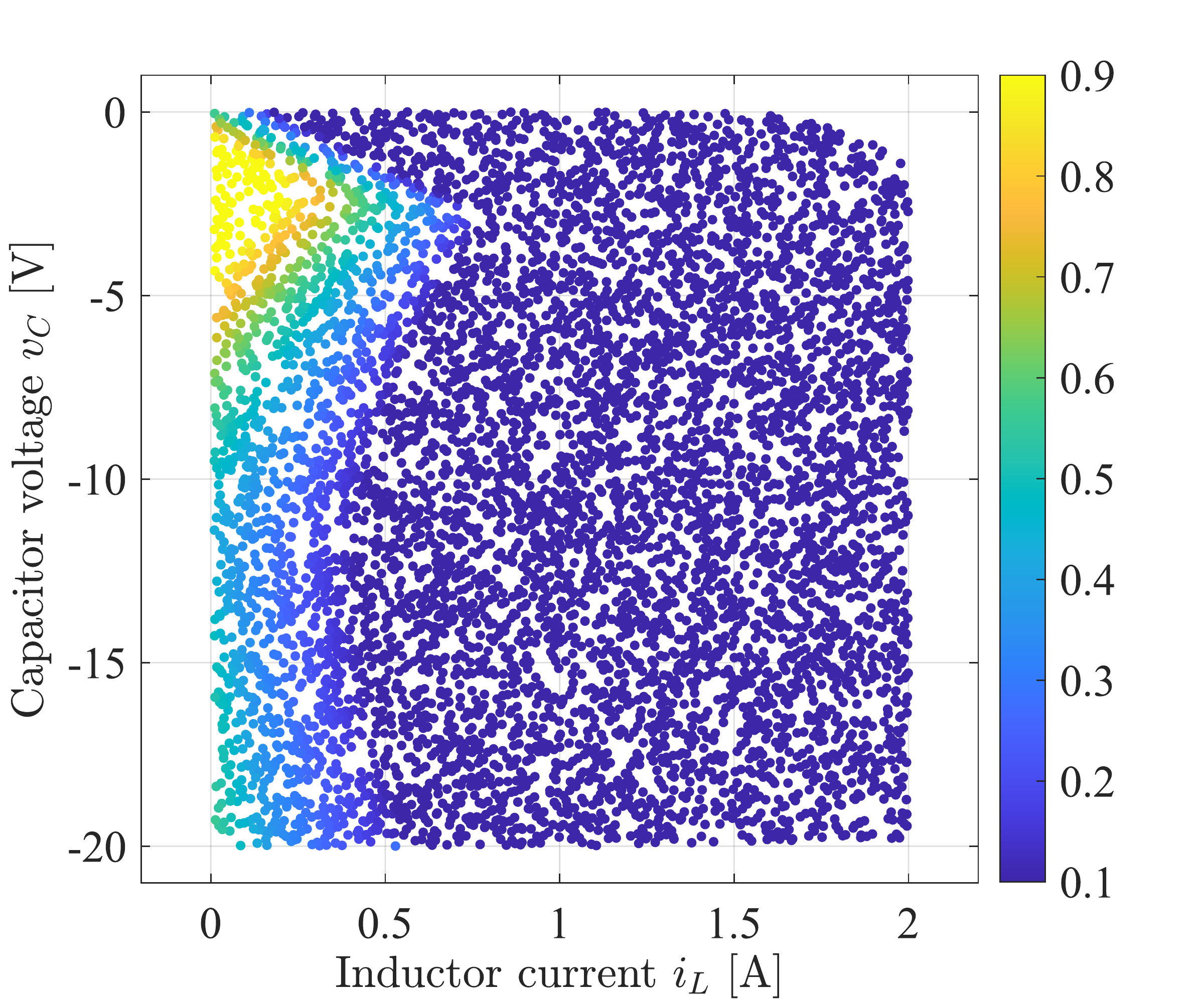}
        \caption{NMPC.}
        \label{fig:1NMPC_Colormap_1}
    \end{subfigure}
    \begin{subfigure}[b]{0.32\linewidth}
        \centering
        \includegraphics[width=\linewidth]{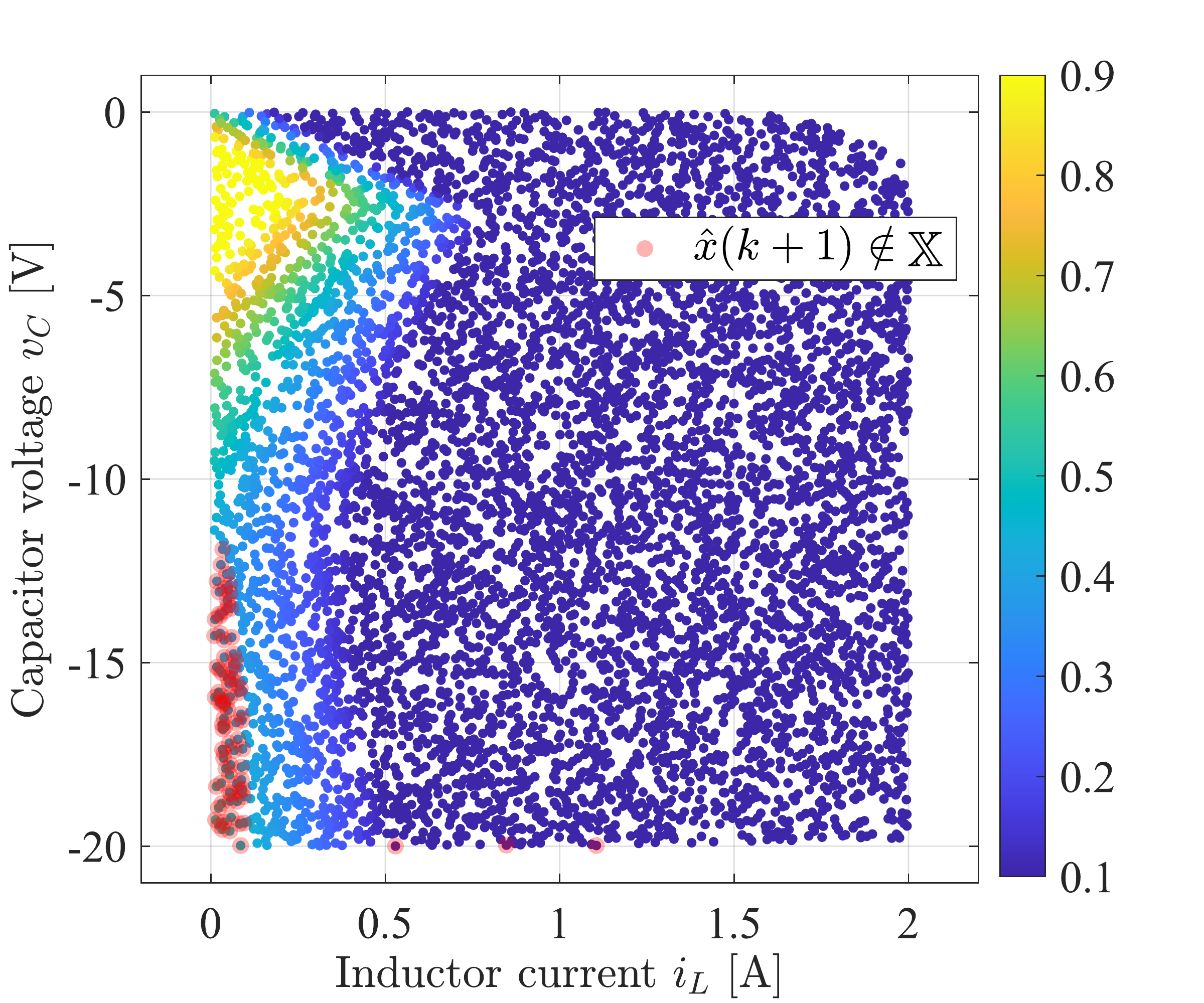}
        \caption{NN-based NMPC.}
        \label{fig:1NMPC_Colormap_2}
    \end{subfigure}
    \begin{subfigure}[b]{0.32\linewidth}
        \centering
        \includegraphics[width=\linewidth]{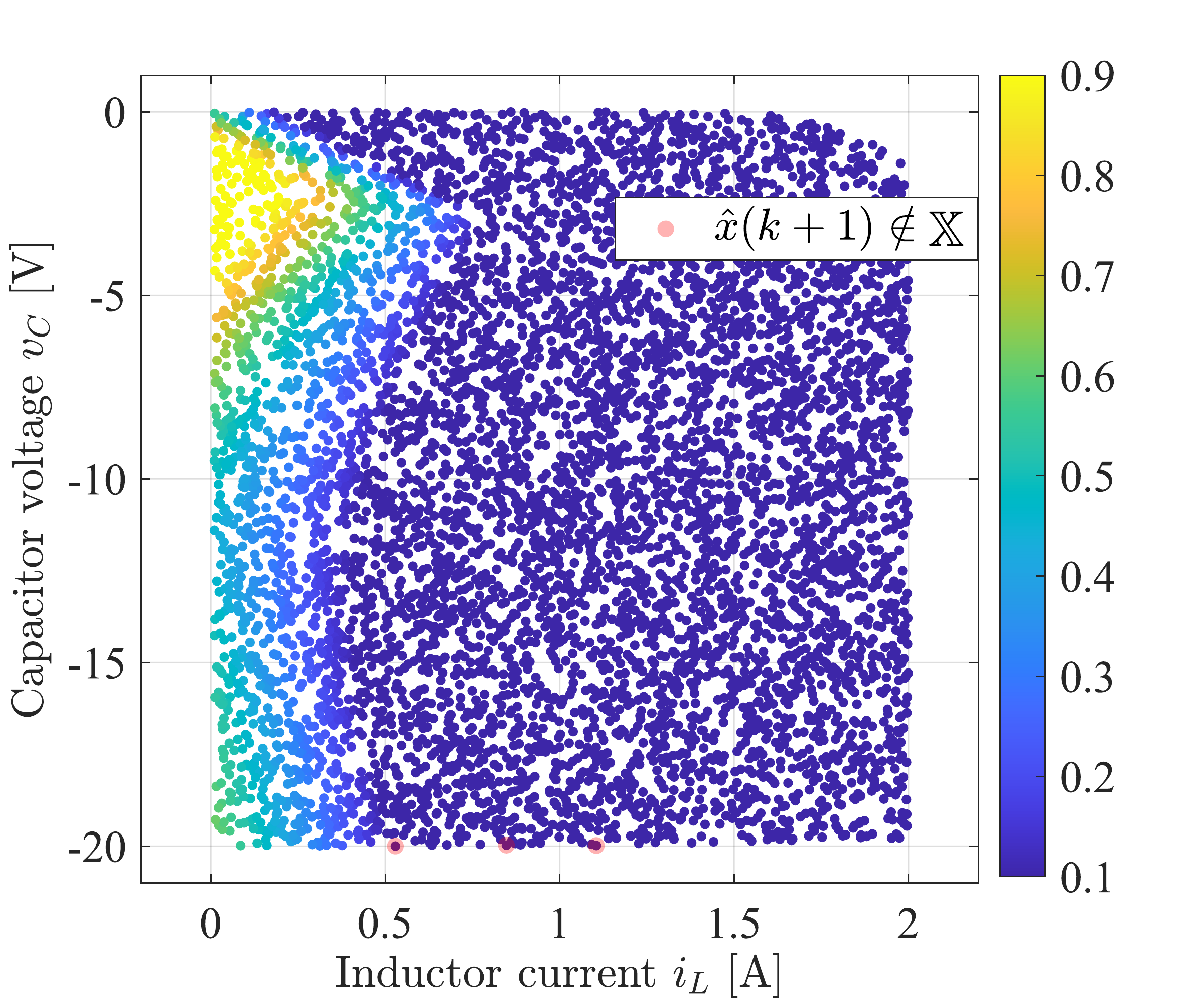}
        \caption{NN-based NMPC ConInf.}
        \label{fig:1NMPC_Colormap_3}
    \end{subfigure}
    \begin{subfigure}[b]{0.32\linewidth}
        \centering
        \includegraphics[width=\linewidth]{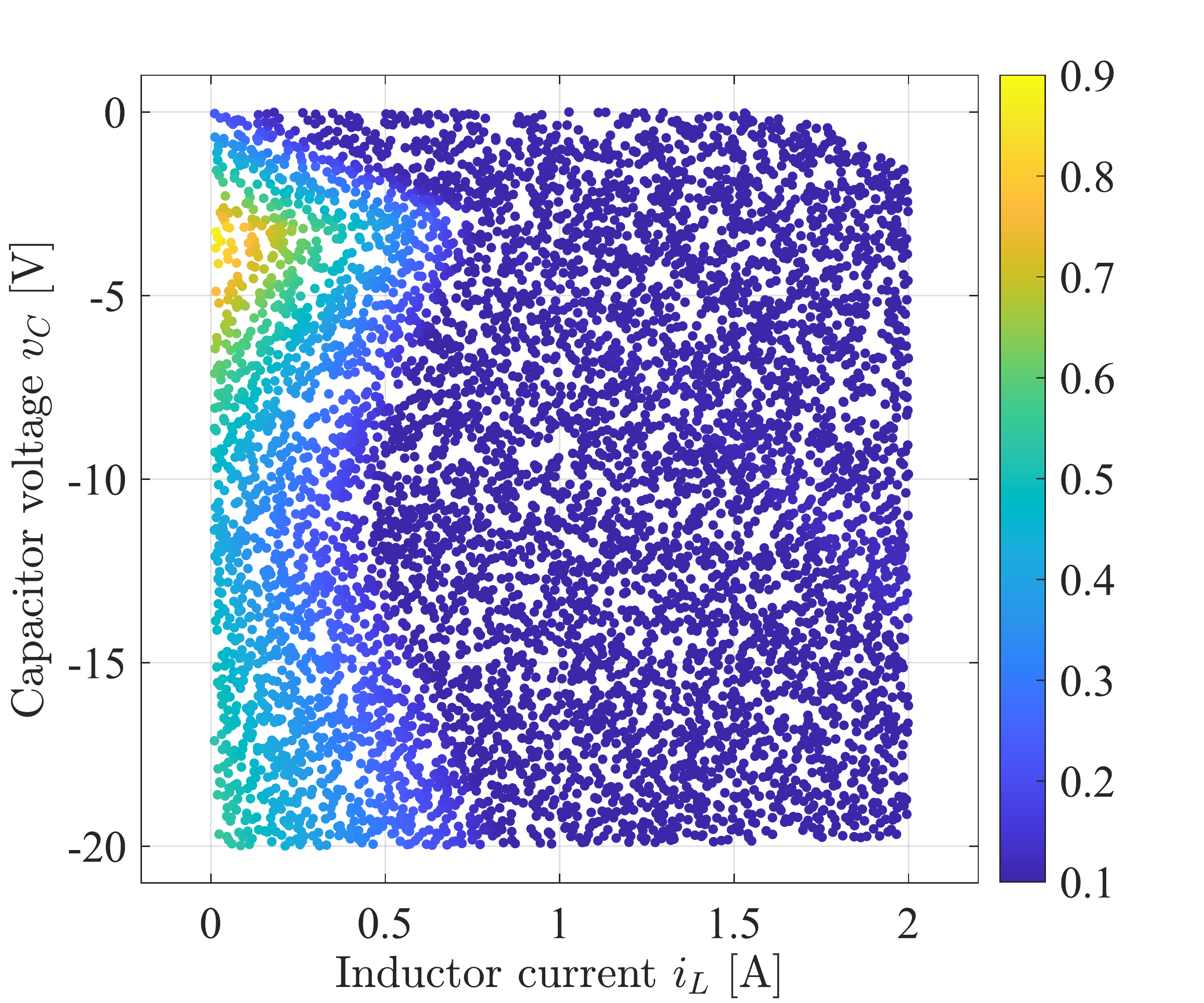}
        \caption{LagNMPC.}
        \label{fig:2LagNMPC_Colormap_1}
    \end{subfigure}
    \begin{subfigure}[b]{0.32\linewidth}
        \centering
        \includegraphics[width=\linewidth]{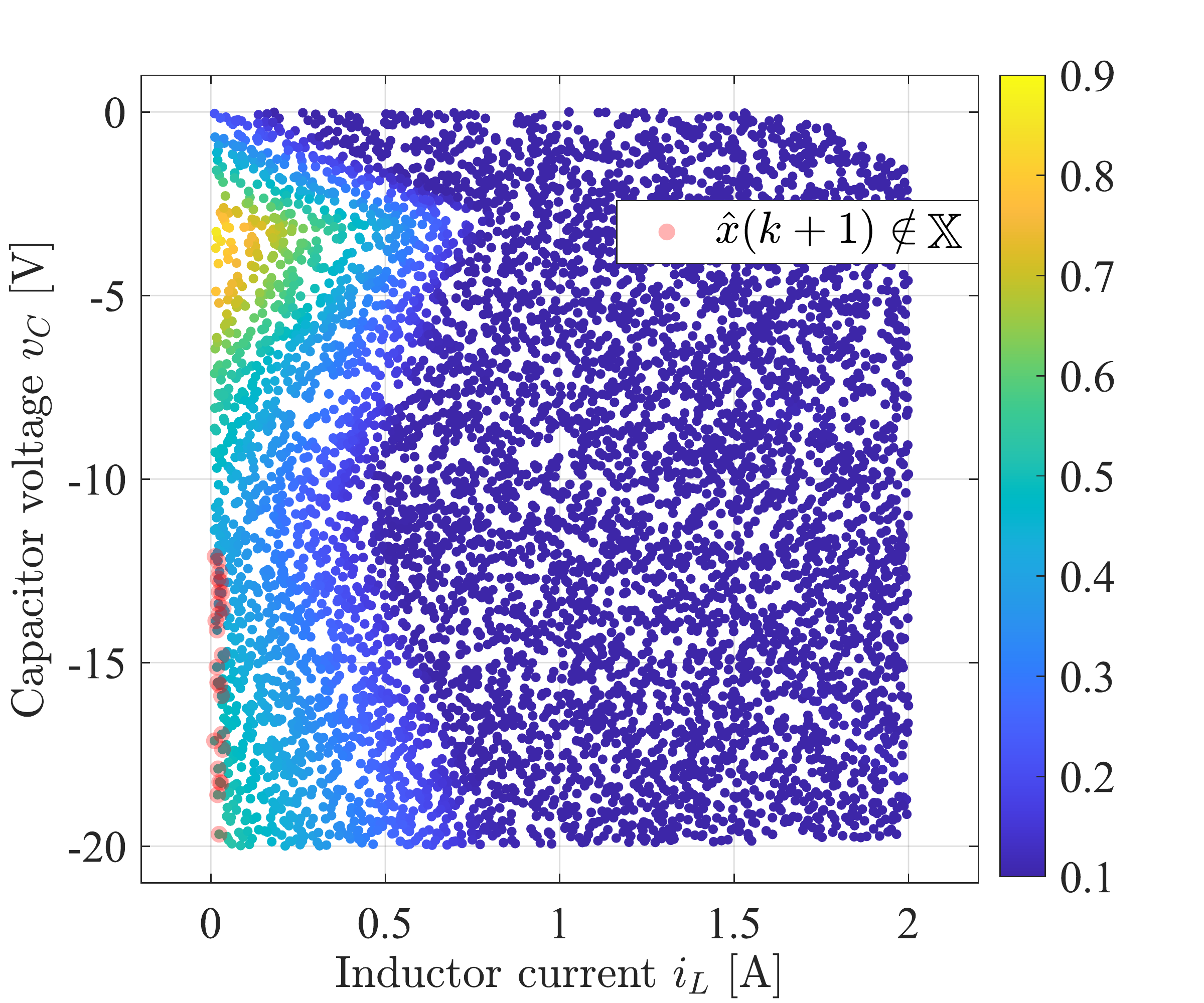}
        \caption{NN-based LagNMPC.}
        \label{fig:2LagNMPC_Colormap_2}
    \end{subfigure}
    \begin{subfigure}[b]{0.32\linewidth}
        \centering
        \includegraphics[width=\linewidth]{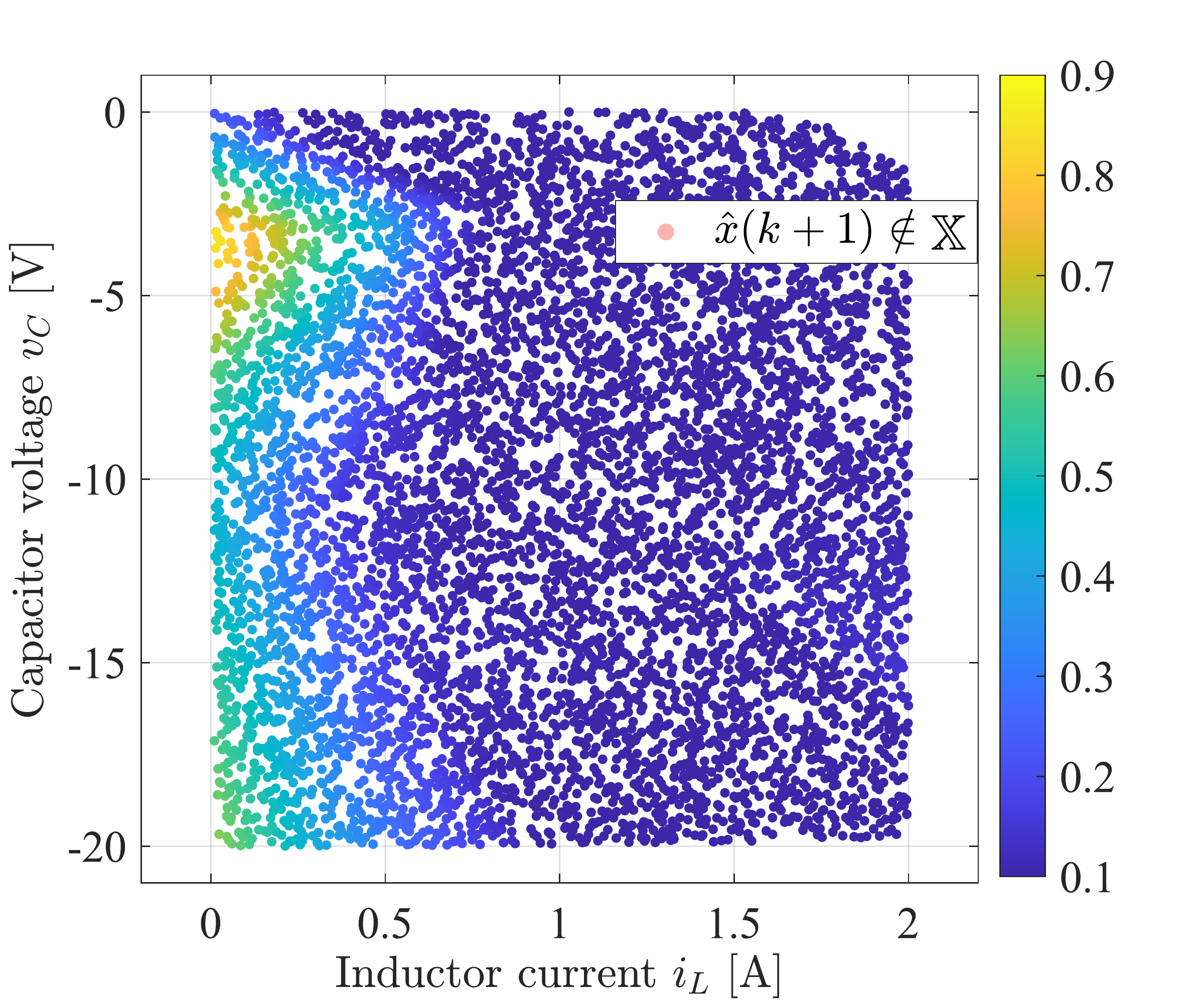}
        \caption{NN-based LagNMPC ConInf.}
        \label{fig:2LagNMPC_Colormap_3}
    \end{subfigure}
    \caption{Comparison of the (Lag)NMPC solutions and the NN-based learnt solutions. The variations of the corresponding inputs over the data set are presented via the scaling color maps. The potential solutions which might violate the state constraint ($x_i^+ \notin \mathbb{X}$) are highlighted with red points.}
    \label{fig:Colormap}
\end{figure*}

Since the considered buck-boost converter has two states and one input, we can intuitively demonstrate the input distribution over the feasible state region as presented in Figure \ref{fig:Colormap}, the variation of the input value is represented via the scaling color map which is restricted between $u_\text{min}$ and $u_\text{max}$.
Figure \ref{fig:1NMPC_Colormap_1} and Figure \ref{fig:2LagNMPC_Colormap_1} illustrate the input action $u_k^*$ performance of the proposed NMPC and LagNMPC, which is calculated using \eqref{eq2:u_NMPC} and \eqref{eq2:u_LagNMPC}. It is observed that by parameterizing the NMPC with Laguerre function, a smoother input distribution is obtained.
Similarly, Figure \ref{fig:1NMPC_Colormap_2} and Figure \ref{fig:2LagNMPC_Colormap_2} illustrate the approximated performance of the corresponding NN-based controllers \eqref{eq3:NNcontroller}. A close color map of the NN-based solution is achieved as compared to that of the (Lag)NMPC solution, 
Still, given the approximated input action, the state at next time instant $x_i^+$ can be predicted. The potential solutions which might violate the state constraint ($x_i^+ \notin \mathbb{X}$) are indicated by red points. As evaluated, the approximation error of the input action would result in violation of the state constraint, which is not expected from the (Lag)NMPC performance.
By introducing the ConInf loss function \eqref{eq3:ConInf} in training the NN, the improved NN-based solutions are shown in Figure \ref{fig:1NMPC_Colormap_3} and Figure \ref{fig:2LagNMPC_Colormap_3}.
As can be seen, all cases that previously resulted in a state constraint violation of $x_1$ have been effectively eliminated.

\begin{figure}[t!]
    \centering
    \begin{subfigure}[b]{1\linewidth}
        \centering
        \includegraphics[width=\linewidth]{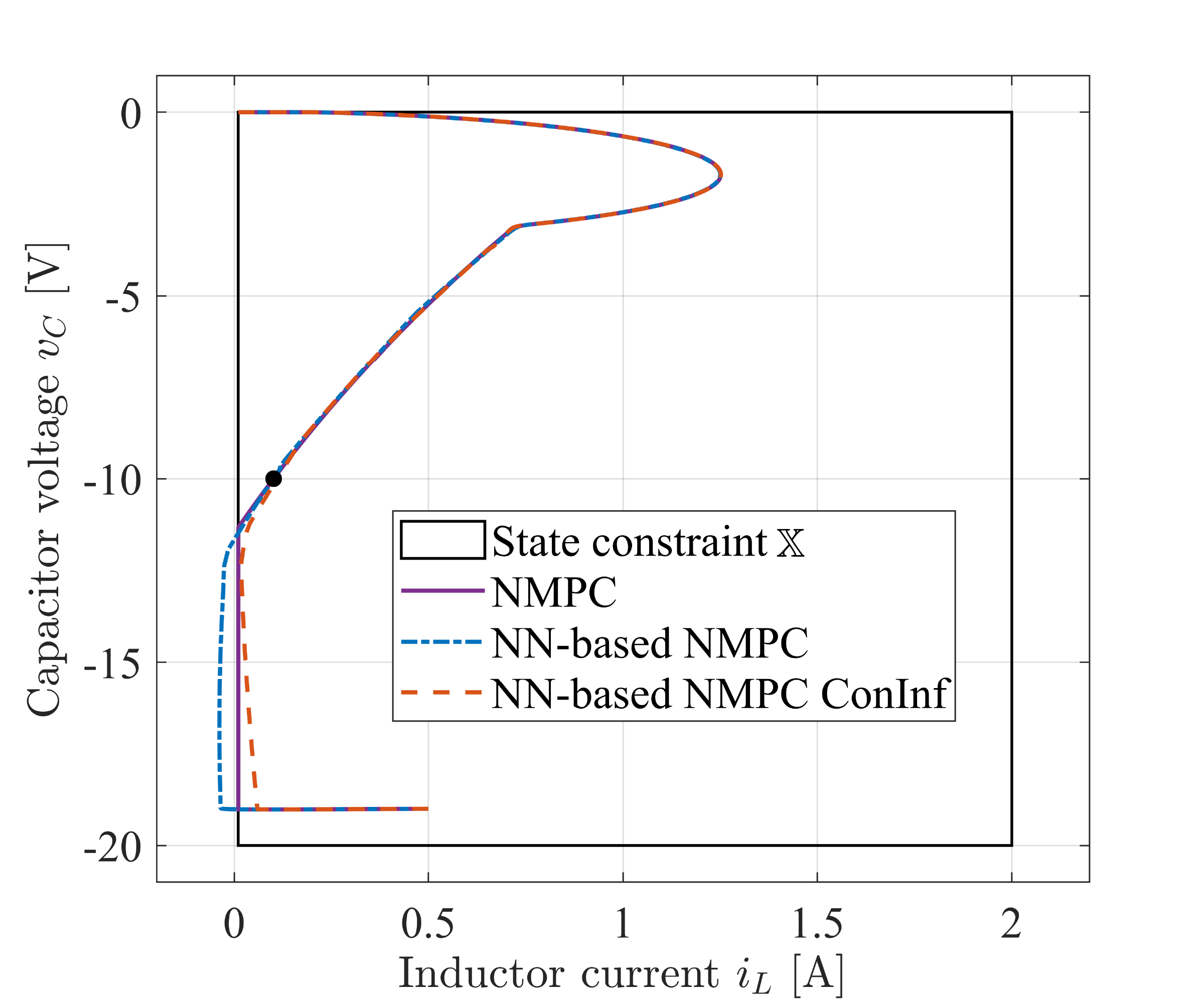}
        \caption{NMPC.}
        \label{fig:1NMPC_x}
    \end{subfigure}
    \begin{subfigure}[b]{1\linewidth}
        \centering
        \includegraphics[width=\linewidth]{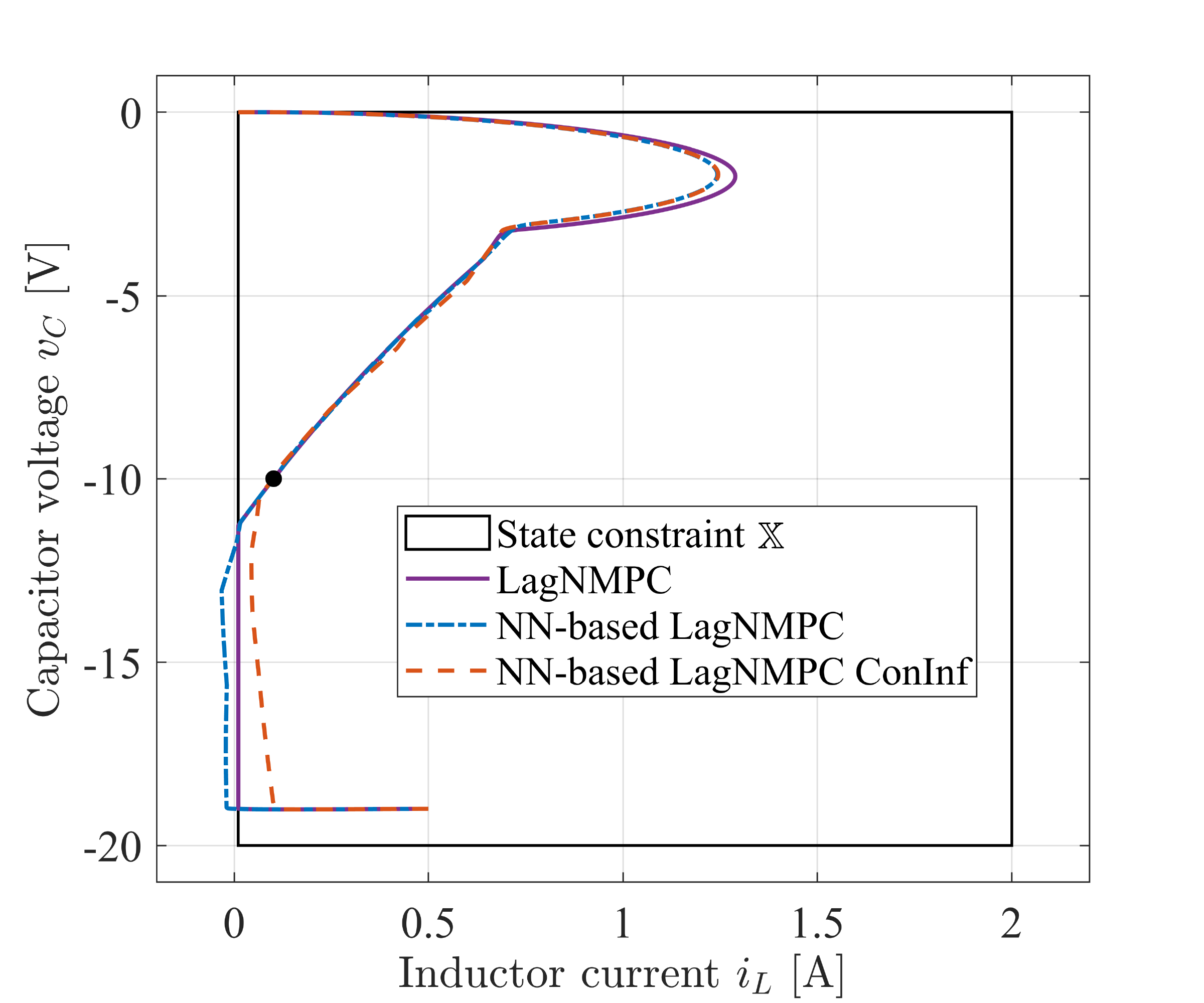}
        \caption{LagNMPC.}
        \label{fig:2LagNMPC_x}
    \end{subfigure}
    \caption{State trajectories with two different initial points.}
    \label{fig:x_traj}
\end{figure}

The closed-loop control performance of all the above-mentioned proposed solutions is also compared.
The offset-free control law is adopted for the NN-based controllers \eqref{eq3:offset_control} with $\epsilon=0.3$.
The evolution of the state trajectories is depicted in Figure \ref{fig:x_traj}. Two initial points ($x_0=[0.01, 0]^\top$ and $x_0=[0.5, -19]^\top$) are selected to demonstrate the control performance. As visualized in Figure \ref{fig:Colormap}, the proposed NN-based control laws would result in state constraint violation performance, see the blue trajectories. Although the constraints-informed NN-based solution does not produce trajectories that closely follows the (Lag)NMPC trajectories, which evolves along the constraint boundary, it ensures that the state moves away from the boundary to prevent violations.

\begin{figure}[t!]
    \centering
    \begin{subfigure}[b]{1\linewidth}
        \centering
        \includegraphics[width=\linewidth]{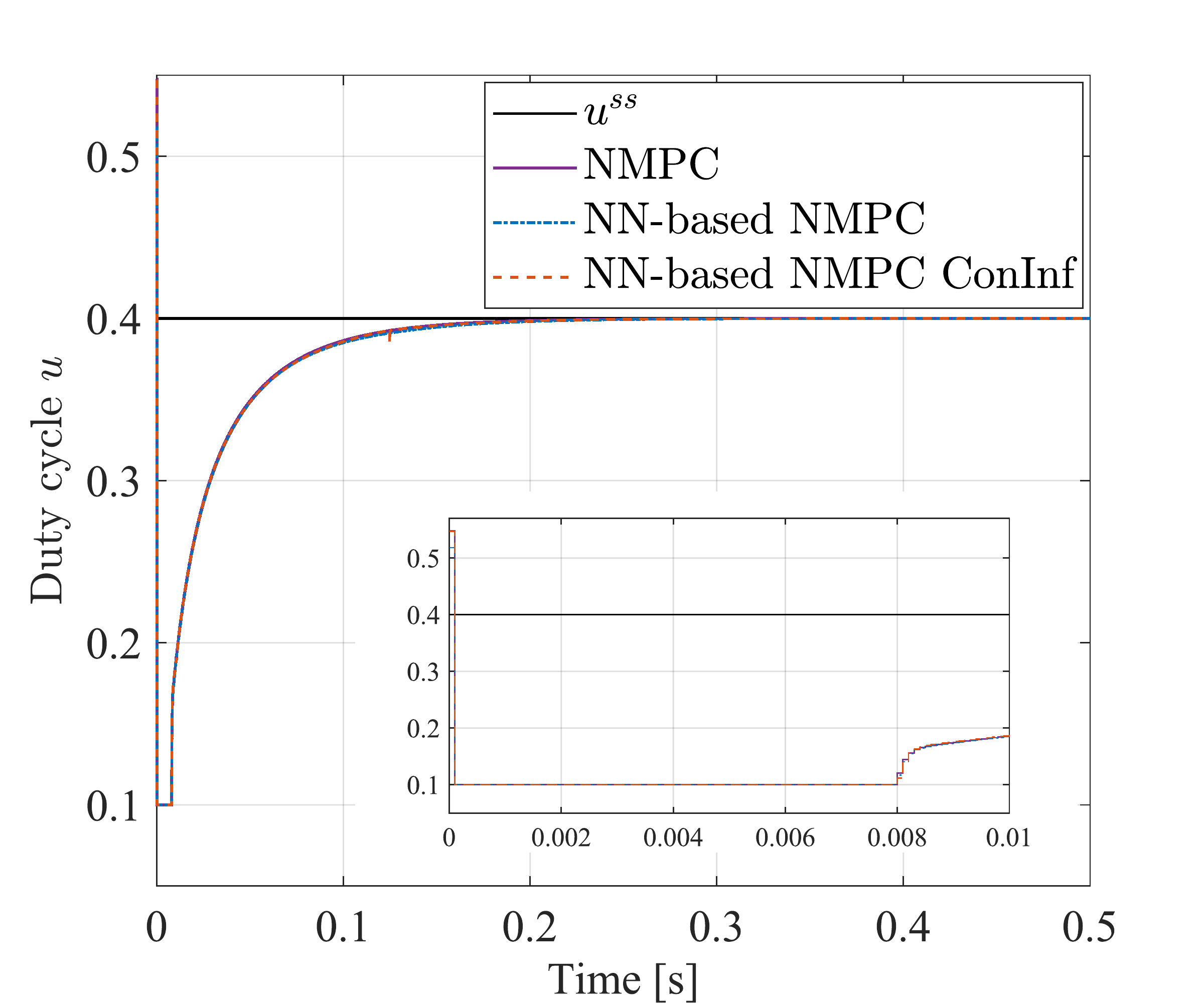}
        \caption{NMPC.}
        \label{fig:1NMPC_u}
    \end{subfigure}
    \begin{subfigure}[b]{1\linewidth}
        \centering
        \includegraphics[width=\linewidth]{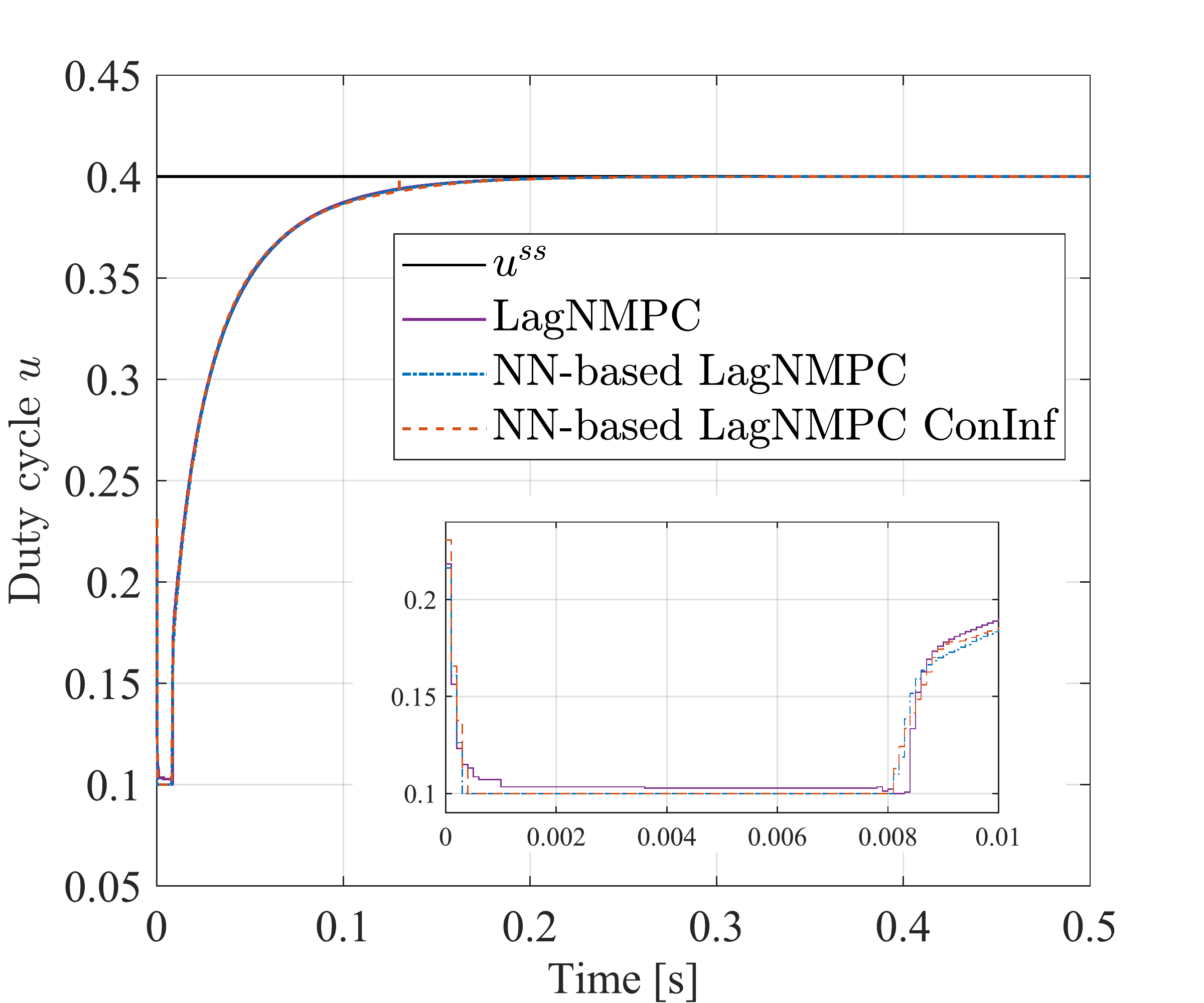}
        \caption{LagNMPC.}
        \label{fig:2LagNMPC_u}
    \end{subfigure}
    \caption{Input time plot with initial point $x_0=[0.01, 0]^\top$.}
    \label{fig:u_time}
\end{figure}

Finally, Figure \ref{fig:u_time} reveals the evolution of all the control inputs, i.e., the duty cycles, over time for the initial condition $x_0=[0.01, 0]^\top$.
Since the proposed NN structures incorporate hard input constraints, the produced results do not violate these limits. 
It is worth mentioning, however, that the spike in the NN-based solution occurs when the offset-free control law becomes active.


In a last step, to obtain insight into the real-time implementation-related aspect of the NN controller, the obtained NN was deployed onto an FPGA architecture. This process followed the high-level synthesis (HLS) procedure, which involves two key steps \cite{chan2021deep}. First, a description of the NN controller was created using a high-level programming language, i.e., C++ in our case. In the second step, this high-level description was converted into hardware description language (HDL) code, which was then used to program the FPGA.

The high-level structure of the network consists of two nested for-loops per layer. These for-loops can be parallelized when converting to HDL code, reducing the computation time. The utilization of hardware resources, however, increases when designing the loops in parallel. The choice of numerics also greatly influences the computation time and resources. Floating-point arithmetic is able to represent a wide range of numbers, due to its dynamic representation. This, however, comes at the cost of accuracy in terms of time and resources. With the use of a fixed-point datatype, throughput time and resources can be decreased in exchange for lower numeric precision. The advantage lies in the fact that this trade-off is part of the design process.
Eventually, fixed-point numerics were used in this work with a $32$-bit word length.
Given the above,  the computation time needed to run the NN controller, synthesized using the Xilinx Vivado HLS tool for Xilinx ZYNQ 7020 SoC with clock frequency of $100\,$MHz, is shown in Table \ref{tab4:computationTimes}. 

\begin{table}[H]
\centering
\caption{FPGA computation time of NN-based controllers.}
\label{tab4:computationTimes}
\begin{tabular}{|c|c|c|}
\hline
Network structure &   No parallel  & Parallel  \\ \hline
NN-based NMPC (Figure \ref{fig2:NN})       &   14.70$\upmu$s      & 2.57$\upmu$s\\
NN-based LagNMPC (Figure \ref{fig:NN_Lag})  & 16.45$\upmu$s &  2.80$\upmu$s \\
\hline
\end{tabular}
\end{table}
\section{Conclusion} \label{Section5}

In this paper, we parameterize the NMPC sequence using Laguerre functions. This results in a smoother control action and renders the potential benefits of approximating the entire control sequence with less variables.
Improved NN structures for approximating LagNMPC with embedded input box constraints are also proposed. Moreover, to mitigate possible state constraint violations caused by the approximation error, a constraints-informed loss function is adopted during the training phase. As a result, the introduction of the aforementioned features makes the proposed approach suitable for applications, such as power electronic systems, where the existence of hard (input) and soft (state) constraints allows for full utilization of the system hardware. 
To demonstrate the effectiveness of the developed solutions, a nonlinear dc-dc buck-boost converter is adopted as a case study. As shown, the proposed approach not only achieves similar performance with online long-horizon NMPC, but it also requires very short execution times, in the range of a few microseconds.


\bibliographystyle{IEEEtran}
\bibliography{SectionBib}

\addtolength{\textheight}{-12cm}   









\end{document}